\documentclass[dvips]{glasnik}
\usepackage{glasnik}
\usepackage{graphicx,amscd}
\usepackage[T1]{fontenc}

\begin{document}

\title[ON PRODUCTS IN THE COARSE SHAPE CATEGORY]{ON PRODUCTS IN THE COARSE SHAPE CATEGORY}
\UniCountry{Ferdowsi University of Mashhad, Iran}
\author[T. Nasri, B. Mashayekhy and H. Mirebrahimi]{Tayyebe Nasri, Behrooz Mashayekhy and Hanieh Mirebrahimi}
\address[Tayyebe Nasri]{Department of Pure Mathematics,\\ Center of Excellence in Analysis on Algebraic Structures,\\ Ferdowsi University of Mashhad,\\
P.O.Box 1159-91775, Mashhad, Iran.}
\email[Tayyebe Nasri]{ta\_na8@stu-mail.um.ac.ir}
\address[Behrooz Mashayekhy]{Department of Pure Mathematics,\\ Center of Excellence in Analysis on Algebraic Structures,\\ Ferdowsi University of Mashhad,\\
P.O.Box 1159-91775, Mashhad, Iran.}
\email[Behrooz Mashayekhy]{bmashf@um.ac.ir}
\address[Hanieh Mirebrahimi]{Department of Pure Mathematics,\\ Center of Excellence in Analysis on Algebraic Structures,\\ Ferdowsi University of Mashhad,\\
P.O.Box 1159-91775, Mashhad, Iran.}
\email[Hanieh Mirebrahimi]{h\_mirebrahimi@um.ac.ir}



\abstract{The paper is devoted to the study of coarse shape of Cartesian products of topological spaces. If the Cartesian product of two spaces $X$ and $Y$ admits an HPol-expansion, which is the Cartesian product of HPol-expansions of these spaces, then $X\times Y$ is a product in the coarse shape category.  As a consequence, the Cartesian product of two compact Hausdorff spaces is a product in the coarse shape category. Finally, we show that the shape groups and the coarse shape groups commute with products under some conditions.
}

\maketitle

\section{Introduction and Preliminaries}
\label{intro}
Since founding shape theory, by K. Borsuk \cite{B0}, it has been developing in several directions. Papers introducing equivalence relations of metric compacta, strictly coarser than the shape equivalence (K. Borsuk q-equivalence, S. Mardesic S-equivalence, S$^*$-equivalence) can be found in the literature since 1976. Recently, N. Kocei\'{c} Bilan and N. Ugle\v{s}i\'{c} \cite{UB} have extended the shape theory by constructing a coarse shape category, denoted by Sh$^*_{(\mathcal{T}, \mathcal{P})}$, where $\mathcal{P}$ is a pro-reflective subcategory of $\mathcal{T}$, whose objects are all objects of $\mathcal{T}$. Its isomorphisms classify objects of $\mathcal{T}$ strictly coarser than the shape does. This category is functorially related to the shape category Sh$_{(\mathcal{T}, \mathcal{P})}$ by a faithful functor and consequently the shape category Sh$_{(\mathcal{T}, \mathcal{P})}$ can be considered as a subcategory of Sh$^*_{(\mathcal{T}, \mathcal{P})}$. Since the homotopy category of polyhedra, HPol, is pro-reflective (dense) in the homotopy category HTop \cite[Theorem 1.4.2]{MS}, the coarse shape category Sh$^*$:=Sh$^*_{(HTop,HPol)}$ is well-defined. The shape and coarse shape coincide on the class of spaces having homotopy type of polyhedra.

Other direction in the developing of the shape theory was a studying of a strong shape. Strong shape for metric compacta was introduced by Edwards and Hastings \cite{EH} and shortly afterward by Bauer \cite{B} for arbitrary spaces. A standard category framework for the strong shape is the strong shape category SSh.

For every category $\mathcal{T}$ an essential question is: Does a category $\mathcal{T}$ admit products?

The product of a pair of objects $X$ and $Y$ (in the categorical sense) is the object $W$ along with the morphisms $p_X : W \rightarrow X$ and $p_Y : W \rightarrow Y$ such that, for every object $Z$ and every pair of morphisms $f : Z \rightarrow X$, $g : Z \rightarrow Y$, there exists a unique morphism $h : Z \rightarrow W$ such that $p_X \circ h = f$ and $p_Y \circ  h= g$.
Although there are several particular results concerning this question for the shape and strong shape categories, this problem is still open for these categories.  It is well-known that for compact Hausdorff spaces $X$ and $Y$, their Cartesian product $X\times Y$ is a product in
the both categories. For ordinary shape this was proved by Keesling \cite{K} and
for strong shape it was proved by Marde\v{s}i\'{c} \cite{M1}. On the other side it was shown
by Keesling that Cartesian product of two non-compact spaces need not be their
product in the shape category. Therefore, the following question naturally has
arisen:\\
(Q) For which topological spaces the Cartesian product $X\times Y$ along
with the induced projections form a product of $X$ and $Y$ in the (strong)
shape category?

In this paper, we consider the problem (Q) in the coarse shape category. In Section \ref{pro}, we intend to verify existence of product in the coarse shape category. Before studying products in coarse shape category we study products in the category pro$^*$-HTop and we show that every pair of inverse systems $\mathbf{X}$, $\mathbf{Y}$ in this category has a product. By a similar argument of \cite{M1} we prove that if the Cartesian product of two spaces $X$ and $Y$ admits an HPol-expansion, which is the Cartesian product of HPol-expansions of these spaces, then $X\times Y$ is a product in the coarse shape category. As a consequence we show that the Cartesian product of two compact Hausdorff spaces is a product in the coarse shape. Also we show that the shape groups and the coarse shape groups (of compact Hausdorff spaces) commute with the product under some conditions.

By the fact that every inverse system is isomorphic to a cofinite inverse system \cite[ Remark 1.1.5]{MS}, in this paper every inverse system is assumed to be cofinite inverse system.

Now, let us recall from \cite{UB} some of the main notions concerning the coarse shape category and a pro$^*$-category. Let $\mathcal{T}$ be a category and let $\mathbf{X}=(X_{\lambda},p_{\lambda\lambda'},\Lambda)$ and $\mathbf{Y}=(Y_{\mu},q_{\mu\mu'},M)$ be two inverse systems in category $\mathcal{T}$. An {\it $S^*$-morphism} of inverse
systems, $(f,f^n_{\mu}): \mathbf{X} \rightarrow \mathbf{Y} $, consists of an index function $f : M \rightarrow\Lambda$ and of a set of $\mathcal{T}$-morphisms $f^n_{\mu}: X_{f(\mu)} \rightarrow Y_{\mu} $, $n\in \Bbb{N}$ and $\mu\in M$ such that, for every related pair $\mu\leq \mu'$ in $M$, there exist a $\lambda\in\Lambda$, $\lambda\geq f(\mu),f(\mu')$, and an $n \in \Bbb{N}$ so that, for every $n'\geq n$, $$q_{\mu\mu'}f_{\mu'}^{n'}p_{f(\mu')\lambda}= f_{\mu}^{n'}p_{f(\mu)\lambda}.$$
If $M=\Lambda$ and $f=1_{\Lambda}$, then $(1_{\lambda},f_{\lambda}^n)$ is said to be a {\it level $S^*$-morphism}.

The {\it composition} of S$^*$-morphisms $(f,f^n_{\mu}): \mathbf{X} \rightarrow \mathbf{Y} $ and $(g,g^n_{\nu}): \mathbf{Y} \rightarrow \mathbf{Z}=(Z_{\nu},r_{\nu\nu'},N) $ is an S$^*$-morphism $ (h,h^n_{\nu})=(g,g^n_{\nu})(f,f^n_{\mu}): \mathbf{X} \rightarrow \mathbf{Z} $, where $h=fg$ and $h^n_{\nu}=g^n_{\nu}f^n_{g(\nu)}$ for all $n\in\Bbb{N}$. The {\it identity S$^*$-morphism} on $\mathbf{X}$ is an S$^*$-morphism $(1_{\Lambda}, 1_{X_{\lambda}}^n): \mathbf{X} \rightarrow \mathbf{X} $, where $1_{\Lambda}$ is the identity function and $1_{X_{\lambda}}^n=1_{X_{\lambda}}$ in $\mathcal{T}$, for all $n\in \Bbb{N}$ and every $\lambda\in\Lambda$.

An S$^*$-morphism $(f,f^n_{\mu}): \mathbf{X} \rightarrow \mathbf{Y} $ is said to be {\it equivalent} to an S$^*$-morphism $(f',f'^n_{\mu}): \mathbf{X} \rightarrow \mathbf{Y} $, denoted by $(f,f^n_{\mu})\sim (f',f'^n_{\mu})$, provided every $\mu\in M$ admits a $\lambda\in\Lambda$ and $n \in \Bbb{N}$ such that $\lambda\geq f(\mu),f'(\mu)$ and for every $n'\geq n$,
$$f_{\mu}^{n'}p_{f(\mu)\lambda}= f'^{n'}_{\mu}p_{f'(\mu)\lambda}.$$

The relation $\sim$ is an equivalence relation among S$^*$-morphisms of inverse systems in $\mathcal{T}$. The {\it category} pro$^*$-$\mathcal{T}$ has as objects all inverse systems $\mathbf{X}$ in $\mathcal{T}$ and as morphisms all equivalence classes $\mathbf{f^*}=[(f,f^n_{\mu})]$ of S$^*$-morphisms $(f,f^n_{\mu})$. The composition in pro$^*$-$\mathcal{T}$ is well defined by putting
\[\mathbf{g^*f^*}=\mathbf{h^*}=[(h,h^n_{\nu})],\]
where $(h,h^n_{\nu})=(g,g^n_{\nu})(f,f^n_{\mu})=(fg, g^n_{\nu}f^n_{g(\nu)})$.
For every inverse system $\mathbf{X}$ in $\mathcal{T}$, the identity morphism in pro$^*$-$\mathcal{T}$ is $\mathbf{1_X^*}=[(1_{\Lambda}, 1^n_{X_{\Lambda}})]$.

In particular if $(X)$ and $(Y)$ are two rudimentary inverse systems in HTop, then every set of mappings $f^n:X\rightarrow Y$, $n\in \Bbb{N}$, induces a map $\mathbf{f^*}:(X)\rightarrow(Y)$ in pro$^*$-HTop.

A functor $\underline{\mathcal{J}}= \underline{\mathcal{J}}_{\mathcal{T}} : pro-\mathcal{T} \rightarrow pro^*-\mathcal{T}$ is defined as follows: For every inverse system $\mathbf{X}$ in $\mathcal{T}$, $\underline{\mathcal{J}} (\mathbf{X}) = \mathbf{X}$ and if $\mathbf{f}\in pro-\mathcal{T}(\mathbf{X}, \mathbf{Y} )$ is represented by $( f , f_{\mu})$, then $\underline{\mathcal{J}} ( \mathbf{f} ) = \mathbf{f^*}=[(f,f^n_{\mu})]\in pro^*-\mathcal{T}(\mathbf{X},\mathbf{Y})$ is represented by the S$^*$-morphism $(f,f^n_{\mu})$, where $f^n_{\mu}= f_{\mu}$ for all $\mu\in M$ and $n\in\Bbb{N}$. Since the functor $\underline{\mathcal{J}}$ is faithful, we may consider the category pro-$\mathcal{T}$ as a subcategory of pro$^*$-$\mathcal{T}$.

Let $\mathcal{P}$ be a subcategory of $\mathcal{T}$. A {\it $\mathcal{P}$- expansion} of an object $X$ in $\mathcal{T}$ is a morphism $\mathbf{p} :X\rightarrow \mathbf{X}$ in pro-$\mathcal{T}$, where $\mathbf{X}$ belongs to pro-$\mathcal{P}$ characterised by the following two properties:\\
(E1) For every object $P$ of $\mathcal{P}$ and every map $h:X\rightarrow P$ in $\mathcal{T}$, there is a $\lambda\in \Lambda$ and a map $f:X_{\lambda}\rightarrow P$ in $\mathcal{P}$ such that $fp_{\lambda}=h$;\\
(E2) If $f_0, f_1:X_{\lambda}\rightarrow P$ in $\mathcal{P}$ satisfy $f_0p_{\lambda}=f_1p_{\lambda}$, then there exists a $\lambda'\geq\lambda$ such that $f_0p_{\lambda\lambda'}=f_1p_{\lambda\lambda'}$.\\
The subcategory $\mathcal{P}$ is said to be {\it pro-reflective} ({\it dense}) subcategory of $\mathcal{T}$ provided every object $X$ in $\mathcal{T}$ admits a $\mathcal{P}$-expansion $\mathbf{p} :X\rightarrow \mathbf{X}$.

Let $\mathcal{P}$ be a pro-reflective subcategory of $\mathcal{T}$. Let $\mathbf{p} :X\rightarrow \mathbf{X}$ and $\mathbf{p'} :X\rightarrow \mathbf{X'}$ be two $\mathcal{P}$-expansions of the same object $X$ in $\mathcal{T}$, and let $\mathbf{q} : Y \rightarrow \mathbf{Y}$ and $\mathbf{q'} : Y \rightarrow \mathbf{Y'}$ be two $\mathcal{P}$-expansions of the same object $Y$ in $\mathcal{T}$. Then there exist two natural (unique) isomorphisms $\mathbf{i}:\mathbf{X}\rightarrow \mathbf{X}'$ and $\mathbf{j}:\mathbf{Y}\rightarrow \mathbf{Y}'$ in pro-$\mathcal{P}$ with respect to $\mathbf{p}$, $\mathbf{p'}$ and $\mathbf{q}$, $\mathbf{q'}$, respectively. Consequently $\underline{\mathcal{J}}(\mathbf{i}):\mathbf{X}\rightarrow \mathbf{X}'$ and $\underline{\mathcal{J}}(\mathbf{j}):\mathbf{Y}\rightarrow \mathbf{Y}'$ are isomorphisms in pro$^*$-$\mathcal{P}$. A morphism $\mathbf{f^*}:\mathbf{X}\rightarrow \mathbf{Y}$ is said to be {\it pro$^*$-$\mathcal{P}$ equivalent} to a morphism $\mathbf{f'^*}:\mathbf{X'}\rightarrow \mathbf{Y'}$, denoted by $\mathbf{f^*}\sim\mathbf{f'^*}$, provided the following diagram in pro$^*$-$\mathcal{P}$ commutes:
\begin{equation}
\label{dia}\begin{CD}
\mathbf{X}@>\underline{\mathcal{J}}(\mathbf{i})>>\mathbf{X'}\\
@VV \mathbf{f^*}V@V \mathbf{f'^*}VV\\
\mathbf{Y}@>\underline{\mathcal{J}}(\mathbf{j})>>\mathbf{Y'}.
\end{CD}\end{equation}

This is an equivalence relation on the appropriate subclass of Mor(pro$^*$-$\mathcal{P}$). Now, the {\it coarse shape category} Sh$^*_{(\mathcal{T},\mathcal{P})}$ for the pair $(\mathcal{T},\mathcal{P})$ is defined as follows: The objects of Sh$^*_{(\mathcal{T},\mathcal{P})}$ are all objects of $\mathcal{T}$. A morphism $F^*:X\rightarrow Y$ is the pro$^*$-$\mathcal{P}$ equivalence class $<\mathbf{f^*}>$ of a mapping $\mathbf{f^*}:\mathbf{X}\rightarrow \mathbf{Y}$ in pro$^*$-$\mathcal{P}$. The {\it composition} of $F^*=<\mathbf{f^*}>:X\rightarrow Y$ and $G^*=<\mathbf{g^*}>:Y\rightarrow Z$ is defined by the representatives, i.e. $G^*F^*=<\mathbf{g^*}\mathbf{f^*}>:X\rightarrow Z$. The {\it identity coarse shape morphism} on an object $X$, $1_X^*:X\rightarrow X$, is the pro$^*$-$\mathcal{P}$ equivalence class $<\mathbf{1_X}^*>$ of the identity morphism $\mathbf{1_X}^*$ in pro$^*$-$\mathcal{P}$.

The faithful functor $\mathcal{J}= \mathcal{J}_{(\mathcal{T},\mathcal{P})} : Sh_{(\mathcal{T},\mathcal{P})}\rightarrow Sh^*_{(\mathcal{T},\mathcal{P})}$ is defined by keeping objects fixed and via the inclusion functor $\underline{\mathcal{J}}= \underline{\mathcal{J}}_{\mathcal{T}} : pro-\mathcal{P} \rightarrow pro^*-\mathcal{P}$.
\begin{remark}\label{Sh}
Let $\mathbf{p}:X\rightarrow \mathbf{X}$ and $\mathbf{q}:Y\rightarrow \mathbf{Y}$ be $\mathcal{P}$-expansions of $X$ and $Y$ respectively. For every morphism $f:X\rightarrow Y$ in $\mathcal{T}$, there is a unique morphism $\mathbf{f}:\mathbf{X}\rightarrow \mathbf{Y}$ in pro-$\mathcal{P}$ such that the following diagram commutes in pro-$\mathcal{P}$:
\begin{equation}
\label{dia}\begin{CD}
\mathbf{X}@<<\mathbf{p}<X\\
@VV \mathbf{f}V@V fVV\\
\mathbf{Y}@<<\mathbf{q}<Y.
\end{CD}\end{equation}
If we take other $\mathcal{P}$-expansions $\mathbf{p'}:X\rightarrow \mathbf{X'}$ and $\mathbf{q'}:Y\rightarrow \mathbf{Y'}$, we obtain another morphism $\mathbf{f'}:\mathbf{X'}\rightarrow \mathbf{Y'}$ in pro-$\mathcal{P}$ such that $\mathbf{f'}\mathbf{p'}=\mathbf{q'}f$ and so we have $\mathbf{f}\sim\mathbf{f'}$ and hence $\underline{\mathcal{J}}(\mathbf{f})\sim \underline{\mathcal{J}}(\mathbf{f'})$ in pro$^*$-$\mathcal{P}$. Therefore every morphism $f\in \mathcal{T}(X,Y)$ yields an pro$^*$-$\mathcal{P}$ equivalence class $<\underline{\mathcal{J}}(\mathbf{f})>$, i.e. a coarse shape morphism $F^*:X\rightarrow Y$, denoted by $\mathcal{S}^*(f)$. If we put $\mathcal{S}^*(X)=X$ for every object $X$ of $\mathcal{T}$, we obtain a functor $\mathcal{S}^*:\mathcal{T}\rightarrow Sh^*$, called the {\it coarse shape functor}.
\end{remark}

Since the homotopy category of polyhedra HPol is pro-reflective (dense) in the homotopy category HTop \cite[Theorem 1.4.2]{MS}, the coarse shape category Sh$^*_{(HTop,HPol)}$=Sh$^*$ is well defined.

\section{Products in coarse shape category}\label{pro}
 As we know the shape category doesn't have the product, in general. Keesling \cite{K} proved that for compact Hausdorff spaces, $X\times Y$ is a product in shape category and Marde\v{s}i\'{c} \cite{M1} showed this result for strong shape category. We intend to verify the existence of product in the coarse shape category. By a similar argument of \cite{M1} we proved that the Cartesian product of two compact Hausdorff spaces is a product in this category. In this section, we prove that if the Cartesian product of two spaces $X$ and $Y$ admits an HPol-expansion, which is the Cartesian product of HPol-expansion of these spaces, then $X\times Y$ is a product in the coarse shape category. Finally, we show that the $k$th shape groups and the $k$th coarse shape groups of compact Hausdorff spaces commute with the product for every $k\in \Bbb{N}$.

Now, we intend to prove  that $X\times Y$ is a product in the coarse shape category under some conditions. First consider some notations.

Let $\mathbf{p}:X\rightarrow \mathbf{X}=(X_{\lambda},p_{\lambda\lambda'},\Lambda)$ and $\mathbf{q}:Y\rightarrow \mathbf{Y}=(Y_{\mu},q_{\mu\mu'},M)$ be morphisms of pro-HTop.
Then $\mathbf{X} \times \mathbf{Y} = (X_{\lambda}\times Y_{\mu},p_{\lambda\lambda'}\times q_{\mu\mu'},\Lambda\times M)$ is an inverse system and $H$-mappings $p_{\lambda}\times q_{\mu}:X \times Y\rightarrow X_{\lambda}\times Y_{\mu}$
form a $S^*$-morphism $\mathbf{p}\times \mathbf{q}:X \times Y\rightarrow \mathbf{X}\times \mathbf{Y}$ in pro-HTop, where $\Lambda\times M$ is directed by putting $(\lambda, \mu)\leq (\lambda' , \mu')$ if and only if $\lambda\leq\lambda'$ and $\mu\leq\mu'$.

To define the canonical projection $\mathbf{\pi^*_X}$ fix an index $\mu\in M$. Let ${}^{\mu}\pi^*_\mathbf{X}: \mathbf{X}\times \mathbf{Y}\rightarrow\mathbf{X}$ be the $S^*$-morphism which consists of the index function $f^{\mu}:\Lambda\rightarrow\Lambda\times M$ defined by $f^{\mu}(\lambda)=(\lambda,\mu)$ and of the set of homotopy classes of projections ${}^{\mu}\pi^{n}_{X_\lambda}:X_{\lambda}\times Y_{\mu}\rightarrow X_{\lambda}$, $\lambda\in \Lambda$, $n\in\Bbb{N}$. If $\mu'$ is another index in $M$, then  ${}^{\mu}\pi^*_\mathbf{X}$ and ${}^{\mu'}\pi^*_\mathbf{X}$ are equivalent $S^*$-morphisms. Hence they induce the same morphism of pro$^*$-HTop which is denoted by $\pi^*_\mathbf{X}$. In the following theorem we show that every pair of inverse systems $\mathbf{X}$ and $\mathbf{Y}$ has a product in the category of pro$^*$-HTop.
\begin{theorem}\label{h}
Let $\mathbf{X}$ and $\mathbf{Y}$ be inverse systems of spaces. Let $\mathbf{f^*}: \mathbf{Z}\rightarrow \mathbf{X} $ and $\mathbf{g^*}: \mathbf{Z}\rightarrow \mathbf{Y}$ be morphisms of pro$^*$-HTop, then there exists a unique morphism  $\mathbf{h^*}: \mathbf{Z}\rightarrow \mathbf{X}\times \mathbf{Y}$ in pro$^*$-HTop such that $\pi^*_\mathbf{X}\mathbf{h^*}=\mathbf{f^*}$ and $\pi^*_\mathbf{Y}\mathbf{h^*}=\mathbf{g^*}$.
\end{theorem}
\begin{proof}
Let $\mathbf{f^*}: \mathbf{Z}=(Z_{\nu},r_{\nu\nu'},N)\rightarrow \mathbf{X}=(X_{\lambda}, p_{\lambda\lambda'}, \Lambda)$ and $\mathbf{g^*}: \mathbf{Z}\rightarrow \mathbf{Y}=(Y_{\mu}, q_{\mu\mu'}, M)$ be given by their representatives $(f,f^n_{\lambda})$ and $(g, g^n_{\mu})$, respectively. To prove existence of $\mathbf{h^*}$, consider the morphism $\mathbf{h^*}: \mathbf{Z}\rightarrow \mathbf{X}\times \mathbf{Y}$ in pro$^*$-HTop represented by $S^*$-morphism $(h, h^n_{(\lambda, \mu)})$ which is defined as follows: Since $f(\lambda), g(\mu)\in N$ and $N$ is directed there exists an $\eta\in N$ such that $\eta\geq f(\lambda), g(\mu)$, we can choose one $\eta$ with this property by the axiom of choice and define the index function $h:\Lambda\times M\rightarrow N$ by $h(\lambda, \mu)=\eta$. For all $(\lambda,\mu)\in \Lambda\times M$ and $n\in\Bbb{N}$ we define $h^n_{(\lambda, \mu)}: Z_{h(\lambda,\mu)}\rightarrow X_{\lambda}\times Y_{\mu}$ by $h^n_{(\lambda, \mu)}(z)=(f^n_{\lambda} r_{f(\lambda)h(\lambda,\mu)}(z), g^n_{\mu}r_{g(\mu)h(\lambda,\mu)}(z))$, where $z\in Z_{\eta}$. We first prove that $\mathbf{h^*}$ is a morphism in pro$^*$-HTop, i.e. $(h, h^n_{(\lambda, \mu)})$ is an $S^*$-morphism. Let $(\lambda, \mu)\leq (\lambda' , \mu')$, then $\lambda\leq\lambda'$ and $\mu\leq\mu'$. Since  $\mathbf{f^*}$ is a morphism in pro$^*$-HTop, for $\lambda\leq\lambda'$ there exist a $\nu_1\in N$, $\nu_1\geq f(\lambda), f(\lambda')$, and an $n_1\in \Bbb{N}$ such that for every $n'\geq n_1$,
\begin{equation}\label{eq1}
f_{\lambda}^{n'}r_{f(\lambda)\nu_1}=p_{\lambda\lambda'}f_{\lambda'}^{n'}r_{f(\lambda')\nu_1}.
\end{equation}
Since $\mathbf{g^*}$ is a morphism in pro$^*$-HTop, for $\mu\leq\mu'$ there exist a $\nu_2\in N$, $\nu_2\geq g(\mu), g(\mu')$, and an $n_2\in \Bbb{N}$ such that for every $n'\geq n_2$,
\begin{equation}\label{eq2}
g_{\mu}^{n'}r_{g(\mu)\nu_2}=q_{\mu\mu'}g_{\mu'}^{n'}r_{g(\mu')\nu_2}.
\end{equation}
Since $N$ is directed there exists a $\nu_3\in N$ such that $\nu_3\geq \nu_1 , \nu_2$. Also there exists a $\nu_4\in N$ such that $\nu_4\geq h(\lambda,\mu), h(\lambda',\mu')$ and so a $\nu\in N$ in which $\nu\geq\nu_3, \nu_4$. By putting $n=max \{n_1, n_2\}$ we have the following equalities for every $n'\geq n$:

\begin{align}
h^{n'}_{(\lambda, \mu)}r_{h(\lambda,\mu)\nu} &=(f^{n'}_{\lambda} r_{f(\lambda)h(\lambda,\mu)}r_{h(\lambda,\mu)\nu}, g^{n'}_{\mu}r_{g(\mu)h(\lambda,\mu)}r_{h(\lambda,\mu)\nu})\nonumber\\
& =(f^{n'}_{\lambda} r_{f(\lambda)\nu}, g^{n'}_{\mu}r_{g(\mu)\nu})\nonumber\\
& =(f^{n'}_{\lambda} r_{f(\lambda)\nu_1}r_{\nu_1\nu}, g^{n'}_{\mu}r_{g(\mu)\nu_2}r_{\nu_2\nu})\nonumber\\
& =(p_{\lambda\lambda'}f_{\lambda'}^{n'}r_{f(\lambda')\nu_1} r_{\nu_1\nu}, q_{\mu\mu'}g_{\mu'}^{n'}r_{g(\mu')\nu_2}r_{\nu_2\nu}) & &\hspace{-2cm} \text{by  (\ref{eq1}) and (\ref{eq2})}, \nonumber\\
& =(p_{\lambda\lambda'}f_{\lambda'}^{n'}r_{f(\lambda')\nu}, q_{\mu\mu'}g_{\mu'}^{n'}r_{g(\mu')\nu})\nonumber\\
& =(p_{\lambda\lambda'}f_{\lambda'}^{n'}r_{f(\lambda')h(\lambda',\mu')} r_{h(\lambda',\mu')\nu}, q_{\mu\mu'}g_{\mu'}^{n'}r_{g(\mu')h(\lambda',\mu')}r_{h(\lambda',\mu')\nu})\nonumber\\
& =(p_{\lambda\lambda'}\times q_{\mu\mu'})h^{n'}_{(\lambda', \mu')}r_{h(\lambda',\mu')\nu}.\nonumber
\end{align}
This proves that $(h, h^n_{(\lambda, \mu)})$ is an $S^*$-morphism. Note that, for $(\lambda,\mu)\in\Lambda\times M$ and for all $n\in \Bbb{N}$, one has ${}^{\mu}\pi_{X_{\lambda}}^{n}(h^n_{(\lambda, \mu)})= f^n_{\lambda}r_{f(\lambda)\nu}$ and hence $\pi^*_{\mathbf{X}}\mathbf{h^*} =\mathbf{f^* }$. Analogously, $\pi^*_{\mathbf{Y}}\mathbf{h^*} = \mathbf{g^*}$.

To prove uniqueness, assume that we have another morphism $\mathbf{h'^*} :\mathbf{Z}\rightarrow \mathbf{X} \times \mathbf{Y}$ is represented by the $S^*$-morphism $(h', h'^n_{(\lambda, \mu)})$ such that
$\pi^*_{\mathbf{X}}\mathbf{h'^*} = \mathbf{f ^*}$ and $\pi^*_{\mathbf{Y}} \mathbf{h'^*} = \mathbf{g^*}$. Since $\pi^*_{\mathbf{X}}\mathbf{h'^*} = \mathbf{f ^*}$ for $(\lambda,\mu)\in\Lambda\times M$ there exist a $\nu_1\in N$, $\nu_1\geq f(\lambda), h'(\lambda,\mu)$, and an $n_1\in \Bbb{N}$ such that for every $n'\geq n_1$,
\begin{equation}\label{eq3}
f_{\lambda}^{n'}r_{f(\lambda)\nu_1}={}^{\mu}\pi_{X_{\lambda}}^{n'}(h'^{n'}_{(\lambda, \mu)})r_{h'(\lambda,\mu)\nu_1}.
\end{equation}
Since $\pi^*_{\mathbf{Y}} \mathbf{h'^*} = \mathbf{g^*}$ for $(\lambda,\mu)\in\Lambda\times M$ there exist a $\nu_2\in N$, $\nu_2\geq g(\mu), h'(\lambda,\mu)$, and an $n_2\in \Bbb{N}$ such that for every $n'\geq n_2$,
\begin{equation}\label{eq4}
g_{\mu}^{n'}r_{g(\mu)\nu_2}={}^{\lambda}\pi_{Y_{\mu}}^{n'}(h'^{n'}_{(\lambda, \mu)})r_{h'(\lambda,\mu)\nu_2}.
\end{equation}
Since $N$ is directed there exists a $\nu_3\in N$ such that $\nu_3\geq \nu_1 , \nu_2$. Also there exists a $\nu \in N$ such that $\nu\geq h(\lambda,\mu),\nu_3$. From (\ref{eq3}) and (\ref{eq4}) we have,
\begin{equation}\label{eq5}
f_{\lambda}^{n'}r_{f(\lambda)\nu}={}^{\mu}\pi_{X_{\lambda}}^{n'}(h'^{n'}_{(\lambda, \mu)})r_{h'(\lambda,\mu)\nu}
\end{equation}
and
\begin{equation}\label{eq6}
g_{\mu}^{n'}r_{g(\mu)\nu}={}^{\lambda}\pi_{Y_{\mu}}^{n'}(h'^{n'}_{(\lambda, \mu)})r_{h'(\lambda,\mu)\nu}.
\end{equation}
By putting $n=max \{n_1, n_2\}$ and using (\ref{eq5}) and (\ref{eq6}) we have the following equalities for every $n'\geq n$:
\begin{align}
h'^{n'}_{(\lambda, \mu)}r_{h'(\lambda,\mu)\nu} &=(f_{\lambda}^{n'}r_{f(\lambda)\nu}, g_{\mu}^{n'}r_{g(\mu)\nu})\nonumber\\
& = h^{n'}_{(\lambda, \mu)}r_{h(\lambda,\mu)\nu}\nonumber.
\end{align}
Hence $\mathbf{h'^*}=\mathbf{h^*}$ and so $\mathbf{h^*}$ is unique.
\end{proof}
\begin{theorem}\label{prod2}
If $X$ and $Y$ admit HPol-expansions $\mathbf{p}:X\rightarrow \mathbf{X}$ and $\mathbf{q}:Y\rightarrow \mathbf{Y}$,  respectively such that $\mathbf{p}\times \mathbf{q}: X\times Y\rightarrow \mathbf{X}\times \mathbf{Y}$ is an HPol-expansion, then $X\times Y$ along with the coarse shape morphisms $\mathcal{S}^*(\pi_X)$ and $\mathcal{S}^*(\pi_Y)$, induced by ordinary projections, is a product in the coarse shape category.
\end{theorem}
\begin{proof}
Let $\pi_X :X\times Y\rightarrow X$ and $\pi_Y :X\times Y\rightarrow Y$ denote the canonical projections. We want to show that $X\times Y$ together with coarse shape
morphisms $\mathcal{S}^*(\pi_X)$ and $\mathcal{S}^*(\pi_Y)$ is a product in Sh$^*$. Let $Z$ be a topological space and let $F^* :Z\rightarrow X$ and $G^*:Z\rightarrow Y$
be coarse shape morphisms. We must prove that there exists a unique coarse shape morphism $H^* :Z\rightarrow X \times Y$ with $\mathcal{S}^*(\pi_X)H^* = F^*$ and $\mathcal{S}^*(\pi_Y )H^* = G^*$.

We will first prove uniqueness of $H^*$. Assume that $H^* :Z\rightarrow X \times Y$ has the desired properties. Consider the HPol-expansions $\mathbf{p}\times \mathbf{q}: X\times Y\rightarrow \mathbf{X}\times \mathbf{Y}$ and $\mathbf{r}:Z\rightarrow \mathbf{Z}$ and let $\mathbf{h^*} :\mathbf{Z}\rightarrow \mathbf{X}\times \mathbf{Y}$ be the morphism of pro$^*$-HTop representing
$H^*$. Similarly, consider the HPol-expansions $\mathbf{p}:X\rightarrow \mathbf{X}$ and $\mathbf{q}:Y\rightarrow \mathbf{Y}$ and let
$\mathbf{f ^*} :\mathbf{Z}\rightarrow \mathbf{X}$ and $\mathbf{g^*} :\mathbf{Z}\rightarrow \mathbf{Y}$ be morphisms of pro$^*$-HTop representing $F^*$ and $G^*$,
respectively. Note that the canonical projection ${}^{\mu}\pi^{n}_{X_\lambda}: X_{\lambda}\times Y_{\mu}\rightarrow X_{\lambda}$ satisfies the equality ${}^{\mu}\pi^{n}_{X_\lambda}(p_{\lambda}\times q_{\mu})=p_{\lambda}\pi_X$, for all $n\in \Bbb{N}$. Therefore $\pi^*_\mathbf{X}(\underline{\mathcal{J}}(\mathbf{p}\times\mathbf{q})) = \underline{\mathcal{J}}(\mathbf{p})\underline{\mathcal{J}}(\pi_X)$, i.e. the following diagram is commutative:
\begin{equation}
\label{dia}\begin{CD}
\mathbf{X}\times \mathbf{Y}@<<\underline{\mathcal{J}}(\mathbf{p}\times \mathbf{q})<X\times Y\\
@VV \pi^*_\mathbf{X} V@V \underline{\mathcal{J}}(\pi_X)VV\\
\mathbf{X}@<<\underline{\mathcal{J}}(\mathbf{p})<X.
\end{CD}\end{equation}

It follows that, by Remark \ref{Sh}, $\mathcal{S}^*(\pi_X)=<\pi^*_\mathbf{X}>$ i.e. the morphism $\pi^*_\mathbf{X}$  represents the coarse shape morphism $\mathcal{S}^*(\pi_X)$ and so the morphism $\pi^*_\mathbf{X}\mathbf{h^*}$ represents the coarse shape morphism $\mathcal{S}^*(\pi_X)H^* = F^*$.
However, we know that $\mathbf{f^* }$ represents $F^*$. Consequently, $\pi^*_\mathbf{X}\mathbf{h^*} = \mathbf{f^*} $.

An analogous argument shows that $\pi^*_\mathbf{Y} \mathbf{h^*} = \mathbf{g^*}$. We now apply Theorem \ref{h} and conclude
that $\mathbf{h^*}$ is unique. This implies that $H^*$ is also unique since it represents $\mathbf{h^*}$.

To prove existence of $H^*$, choose $\mathbf{f^*}$ and $\mathbf{g^*}$ as above and let $\mathbf{h^*}$ be the morphism defined in Theorem \ref{h}. Then define
$H^* :Z\rightarrow X \times Y$ as the coarse shape morphism which is represented
by $\mathbf{h^*}$. Arguing as before, $\pi^*_\mathbf{X}\mathbf{h^*}$ is
associated with $\mathcal{S}^*(\pi_X)H^*$. By Theorem \ref{h}, $\pi^*_\mathbf{X}\mathbf{h^*} = \mathbf{f^*}$. Consequently,
$\mathcal{S}^*(\pi_X)H^*$ is associated with $\mathbf{f^*}$. Since $F^*$ is also associated with $\mathbf{f^* }$, we conclude that
$\mathcal{S}^*(\pi_X)H^* = F^*$. Analogously, $\mathcal{S}^*(\pi_Y )H^* = G^*$.
\end{proof}

Marde\v{s}i\'{c} showed that if $\mathbf{p}:X\rightarrow \mathbf{X}$ and $\mathbf{q}:Y\rightarrow \mathbf{Y}$ are HPol-expansions of compact Hausdorff spaces $X$ and $Y$, respectively, then $\mathbf{p}\times \mathbf{q}: X\times Y\rightarrow \mathbf{X}\times \mathbf{Y}$ is also an HPol-expansion of $X\times Y$ \cite[Lemma 2 and Theorem 4]{MS}. Therefore we have the following result from Theorem \ref{prod2}.
\begin{corollary}\label{pro3}
If $X$ and $Y$ are compact Hausdorff spaces, then $X\times Y$ together with the coarse shape morphisms $\mathcal{S}^*(\pi_X)$ and $\mathcal{S}^*(\pi_Y)$ is a product in the coarse shape category Sh$^*$.
\end{corollary}

N. Kocei\'{c} Bilan \cite{B1} introduced the $k$th coarse shape group $\check{\pi}^*_k(X,x)$, $k\in \Bbb{N}$,  as the set of all coarse shape morphisms $F^*:(S^k,*)\rightarrow (X,x)$ with the following multiplication which makes it a group.
$$ F^*+G^*= <\mathbf{f^*}>+<\mathbf{g^*}>=<\mathbf{f^*}+\mathbf{g^*}>=<[(f^n_{\lambda})]+[(g^n_{\lambda})]>=<[(f_{\lambda}^n+g_{\lambda}^n)]>,$$
where coarse shape morphisms $F^*$ and $G^*$ are represented by morphisms $\mathbf{f^*}=[(f, f^n_{\lambda})]$ and $\mathbf{g^*}=[(g, g^n_{\lambda})]:(S^k,*)\rightarrow (\mathbf{X},\mathbf{x})$ in pro$^*$-HPol$_*$, respectively. In follow, we show that the shape groups and the coarse shape groups of coarse shape path connected, compact and Hausdorff spaces commute with the product. Note that N. Kocei\'{c} Bilan proved that if $X$ is a coarse shape path connected space, then  $\check{\pi}^*_k(X,x_0)\cong\check{\pi}^*_k(X,x_1)$ for any two points $x_0, x_1\in X$  and every $k\in \Bbb{N}_0$ \cite[Corollary 1]{B2}.
  \begin{theorem}\label{product}
 Let $X$ and $Y$ be two coarse shape path connected spaces. If $X$ and $Y$ admit HPol-expansions $\mathbf{p}:X\rightarrow \mathbf{X}$ and $\mathbf{q}:Y\rightarrow \mathbf{Y}$, respectively such that $\mathbf{p}\times \mathbf{q}: X\times Y\rightarrow \mathbf{X}\times \mathbf{Y}$ is an HPol-expansion, then $\check{\pi}_k^*(X\times Y)\cong\check{\pi}_k^*(X)\times \check{\pi}_k^*(Y)$, for every $k\in \Bbb{N}$.
\end{theorem}
\begin{proof}
Let $\mathcal{S}^*(\pi_X):X\times Y\rightarrow X$ and $\mathcal{S}^*(\pi_Y):X\times Y\rightarrow Y$ be the induced coarse shape morphisms of canonical projections and assume that $\phi_X:\check{\pi}^*_k(X\times Y)\rightarrow \check{\pi}^*_k(X)$ and $\phi_Y:\check{\pi}^*_k(X\times Y)\rightarrow \check{\pi}^*_k(Y)$ are the induced homomorphisms of $\mathcal{S}^*(\pi_X)$ and $\mathcal{S}^*(\pi_Y)$, respectively. Then there is a homomorphism $\phi=(\phi_X, \phi_Y):\check{\pi}^*_k(X\times Y)\rightarrow \check{\pi}^*_k(X)\times\check{\pi}^*_k(Y)$. By Theorem \ref{prod2}, $X\times Y$ is a product in Sh$^*$ hence we can define a homomorphism $\psi:\check{\pi}^*_k(X)\times\check{\pi}^*_k(Y)\rightarrow\check{\pi}^*_k(X\times Y)$ by $\psi(F^*, G^*)=\lfloor F^*, G^*\rfloor$, where $\lfloor F^*, G^*\rfloor:S^k\rightarrow X\times Y$ is the unique coarse shape morphism such that $\mathcal{S}^*(\pi_X)\lfloor F^*, G^*\rfloor=F^*$ and $\mathcal{S}^*(\pi_Y)\lfloor F^*, G^*\rfloor=G^*$. Indeed, if $F^*=<\mathbf{f^*}=[(f, f^n_{\lambda})]>$ and $G^*=<\mathbf{g^*}=[(g, g^n_{\mu})]>$, then $\lfloor F^*, G^*\rfloor=<\lfloor \mathbf{f^*}, \mathbf{g^*}\rfloor>$, where $\lfloor\mathbf{f^*}, \mathbf{g^*}\rfloor$ is given by $\lfloor f, g\rfloor^n_{\lambda\mu}=f^n_{\lambda}\times g^n_{\mu}: S^k\rightarrow X_{\lambda}\times Y_{\mu}$ (see Theorem \ref{prod2}). The homomorphism $\psi$ is well-defined. Suppose  $F^*$, $F'^*\in \check{\pi}_k^*(X)$ represented by $\mathbf{f^*}$ and $\mathbf{f'^*}$, respectively and $G^*$, $G'^*\in \check{\pi}_k^*(Y)$ associated by $\mathbf{g^*}$ and $\mathbf{g'^*}$, respectively such that $F^*=F'^*$ and  $G^*=G'^*$, i.e. $\mathbf{f^*}\sim \mathbf{f'^*}$ and $\mathbf{g^*}\sim \mathbf{g'^*}$. One can see that $\lfloor\mathbf{f^*}, \mathbf{g^*}\rfloor\sim \lfloor\mathbf{f'^*}, \mathbf{g'^*}\rfloor$ and hence $\lfloor F^*, G^*\rfloor=\lfloor F'^*, G'^*\rfloor$. It is routine to check that $\phi\circ\psi=id$ and $\psi\circ\phi=id$.
\end{proof}
\begin{remark}
(i) It is known that, for a shape path connected space, shape homotopy groups do not depend on the choice of a base point
(see \cite{un}). With similar to proof of Theorem \ref{product} one can show that if $X$ and $Y$ are two shape path connected spaces and if $X$ and $Y$ admit HPol-expansions $\mathbf{p}:X\rightarrow \mathbf{X}$ and $\mathbf{q}:Y\rightarrow \mathbf{Y}$, respectively such that $\mathbf{p}\times \mathbf{q}: X\times Y\rightarrow \mathbf{X}\times \mathbf{Y}$ is an HPol-expansion, then $\check{\pi}_k(X\times Y)\cong\check{\pi}_k(X)\times \check{\pi}_k(Y)$, for every $k\in \Bbb{N}$.\\
(ii) If $X$ and $Y$ are coarse shape path connected, compact and Hausdorff spaces, then by a result of \cite[Theorem 10]{M1} $\mathbf{p}\times \mathbf{q}: X\times Y\rightarrow \mathbf{X}\times \mathbf{Y}$ is an HPol-expansion, where $\mathbf{p}:X\rightarrow \mathbf{X}$ and $\mathbf{q}:Y\rightarrow \mathbf{Y}$ are HPol-expansions of $X$ and $Y$, respectively. Therefore in this case the $\check{\pi}_k$ and $\check{\pi}_k^*$ commute with finite products, for every $k\in \Bbb{N}$, by Theorem \ref{product}.
\end{remark}

\end{document}